\documentclass[10pt]{amsart} %\swapnumbers 
\begin{document}

\title{Uniform Growth of Polycyclic Groups}
\author{Roger C. Alperin}
\address{Department of Mathematics and Computer Science,\ San Jose State
University,\  San Jose, CA 95192 USA } \email{alperin@mathcs.sjsu.edu}

\newtheorem{theorem}{Theorem}[section]
\newtheorem{proposition}[theorem]{Proposition}
\newtheorem{lemma}[theorem]{Lemma} 
\newtheorem{corollary}[theorem]{Corollary}
\newtheorem{remark}[theorem]{Remark}
\newtheorem{definition}[section]{Definition}

\def\proof{{\bf {\noindent}Proof: }} \newsymbol\bsq 1004 
%\newsymbol\bsq\square

\def\endproof{\hfill{$\bsq$}\par\bigskip}
\def\sdirect{\times\hspace{-.18em}\rule{.22pt}{1.05ex}} \def\bar{\overline}

\newcommand\Z{{\bf Z}} 
\newcommand\Q{{\bf Q}} 
\newcommand\C{{\bf C}}
\newcommand\R{{\bf R}}

\def\A{\mathcal A} 
\def\B{\mathcal B} 
\def\K{\mathcal K} 
\def\H{\mathcal H}
\def\L{\mathcal L} 
\def\M{\mathcal M} 
\def\N{\mathcal N}
\def\V{\mathcal V} 
\def\U{\mathcal U}

\def\c{\kappa} 
\def\d{\delta} 
\def\l{\lambda} 
\def\e{\epsilon} 
\def\G{\Gamma}
\def\g{\gamma}

%\keywords{uniform exponential growth, polycyclic}
%\classification{20E99, 20F16, 20F65}

\maketitle 

\centerline{To John Stallings on the occassion of his $65^{th}$ birthday}

\begin{abstract} We prove that polycyclic groups have uniform exponential or
polynomial growth.
\end{abstract}

\section{Introduction}

The Milnor-Wolf Theorem characterizes the finitely generated solvable groups
which have exponential growth; a finitely generated solvable group has
exponential growth iff it is not virtually nilpotent. In fact, Wolf, \cite{W},
showed that a finitely generated nilpotent by finite group has polynomial
growth; then extended this by proving that polycyclic groups which  are not
virtually nilpotent have exponential growth. On the other hand, Milnor,
\cite{M}, showed that finitely generated solvable groups which are not
polycyclic have exponential growth. In both approaches exponential growth can
be deduced from the existence of a free semigroup, \cite{A, R}.

In this article we elaborate on these results by proving that the growth rate
of a polycyclic group $\G$ of exponential growth is uniformly exponential.
This means that the base of the rate of exponential growth $\beta(S,\G)$ is
bounded away from 1, independent of the set of generators, $S$; that is, there
is a constant $\beta(\G)$ so that $\beta(S,\G)\ge \beta(\G)>1$ for any finite
generating set. The growth rate is also related to the spectral radius
$\mu(S,G)$ of the random walk on the Cayley graph, with the given set of
generators, \cite{Gk}.

The exponential polycyclic groups are an important class of groups for
resolving the question of whether or not exponential growth is the same as
uniform exponential growth, since they are in a sense very close to groups of
polynomial growth. The ideas used here for polycyclic groups take advantage of
their linear and arithmetic properties. These may be important tools for
proving uniform growth for other classes of groups. Other methods for proving
uniform growth take advantage of special properties of presentations; for
example, an excess of the number of generators over relations by at least 2
ensures the existence of a subgroup of finite index which maps onto a
non-abelian free group. Here, of course we can not map to free groups, but we
are in  a sense able to map to a non-abelian free semigroup.

\section{Exponential Growth} The growth of a group is measured in the
following way. Choose a finite generating set, $S$, for the group $\G$; define
the length of an element as $\lambda_S(g)=min\{\ n\ |\ g=s_1\cdots s_n, s_i\in
S\cup S^{-1}\}$. The growth function $\beta_n(\G,S)=|\{ g\ |\  \lambda_S(g)\le
n\}|$ depends on the chosen generating set, but as Milnor showed, whether or
not it is exponential is independent of the generating set. A group has
exponential growth if the growth rate, $\beta(S,\G)=limit_{n\rightarrow
\infty} \beta_n(S,\G)^{\frac{1}{n}}$ is strictly greater than 1, for some
(any) generating set. The growth rate is the reciprocal of the radius of
convergence of the growth series, the generating function with coefficients
given by $\beta_n(S,\G)$.

If we have another finite generating set $T=\{t_1,\cdots, t_m\}$ for $\G$, and
$max_j\lambda_S(t_j)\le L$, $max_i\lambda_T(s_i) \le L$, then
$\beta_n(S,\G)\le \beta_{Ln}(T,\G)$ and also the symmetric inequality. It then
follows that $\beta(S,\G)^L\le\beta(T,\G)$ and $\beta(T,\G)^L\le\beta(S,\G)$
so that exponential rate for one set of generators means exponential rate for
any system of generators.

Suppose that $\G$ has a free semigroup on two generators, say $a,b$. Let
$T=S\cup \{a,b\}$; then, we have the interesting feature, using similar estimates as
above, that $\beta(S,\G)^L\ge \beta(T,\G)\ge 2$, where $L$ is the maximum of
the length of the elements $a, b$ in terms of the generators $S$.

For a group with exponential growth we consider  $$\beta(\G)=inf_{S}
\beta(S,\G).$$ If $\beta(\G)>1$ then $\G$ is said to have uniform
(exponential) growth.

\begin{proposition} \label{bound} Given a finitely generated group $\G$.
Suppose there is a constant $L$ so that for every finite generating set $S$
there are two elemnts $a_S, b_S$ which generate a free semigroup and whose
lengths in that generating set are  bounded by $L$, then $\G$ has uniform
exponential growth. \end{proposition}

\proof It follows immediately from above,  $\beta(S,\G)\ge 2^{\frac{1}{L}}$.
\endproof

The following uses the fact that given a generating set $S$ for $\G$ the set
of elements of the subgroup $\K$, a subgroup of index $d$, which are words in
$S$ of length at most $2d-1$ give a generating set for $\K$.

\begin{proposition}[\cite{SW}] If a group $\G$ has a subgroup $\K$ of finite
index $d$, then $\beta(\G)\ge \beta(\K)^{\frac{1}{2d-1}}$. \label{SW}
\end{proposition}

The following is immediate from the definition and Proposition \ref{SW} on
passage to finite index. \begin{corollary} If a group $\G$ has a subgroup of
finite index which has a homomorphic image which has  uniform growth, then
$\G$ has uniform growth. \label{quo} \end{corollary}

Let $A$ be an automorphism of infinite order of a finitely generated abelian
group $\V$. The group $\V_A=\V$ $\sdirect$ $_A\Z$ is the split extension of
$\V$ by the automorphism $A$. In the split extension, $\G_A$, let $t$ denote 
the element  which maps to the generator of $\Z$, and which acts via
conjugation on $\V$ via the automorphism $A$. Let $M(v)$ denote the cyclic
submodule generated by $v\in\V$. Let $\V\otimes_{\Z}{\C}$ denote the extension
of scalars to the complex numbers, $\C$.

\begin{lemma} \label{free}Let $A$ be an automorphism of a finitely generated
abelian group $\V$.  If $A$ has an eigenvalue $|\l|\ge 2$ on
$M(v)\otimes_{\Z}{\C}$ then the split extension, $\V_A$, contains the free
semigroup generated by $tv$ and $v$. \end{lemma} \proof Let $v_i=t^ivt^{-i}$.
Consider the set $W=\{w\ |\ w=\prod_{i=1}^n {v_{i}^{e_i}}, e_i\in\{ 0, 1\} ,
n\ge 1 \}$.  Let $U=\{t^i\ |\ i\ge 1\}$.  Set $X=WU$. Consider the sets $tvX$
and $tX$. By a calculation,  $t(wt^i)=t(\prod_{i=1}^n
{v_{i}^{e_i}})t^{-1}t^{i+1}=\prod_{i=1}^{n} {v_{i+1}^{e_i}}t^{i+1}$  and also
$tv(wt^i)=tvt^{-1}twt^{-1}t^{i+1}=v_1\prod_{i=1}^{n} {v_{i+1}^{e_i}}t^{i+1}$;
thus, $tX\cup tvX\subseteq X$.

Furthermore, these two subsets $tvX$ and $tX$ do not meet. For if they did
then there would be an equality $tw_1t^i=tvw_2t^j$ for some $w_1, w_2\in W$.
From the calculations above and the uniqueness of the power of $t$ in $U$ it
follows that $i=j$ and then also that $tw_1t^{-1}=v_1(tw_2t^{-1})$. This last
equation gives rise to two polynomials $W_1, W_2$ in $A$, the matrix of
conjugation by $t$, such that $AW_1(A)v=Av+AW_2(A)v$. But since $v$ generates
the cyclic module this means $AW_1(A)-AW_2(A)=A$ or simply $W_1(A)-W_2(A)=I$
where the coefficients of $W_1, W_2$ are 0 or 1, without constant term.
Consequently this identity is also valid for every eigenvalue $\lambda$ of $A$
and therefore $W_1(\l)-W_2(\l)=1$ so that for some $m\ge 1$ (highest power in
$W_1$ or $W_2$), $\l^m=\e_0+\e_1\l^1+\cdots+\e_{m-1}\l^{m-1}$, $\e_j\in \{0,
\pm1\}$, so that $$|\l|^m\le 1+|\l|+|\l|^2+\cdots+|\l|^{m-1}\le
\frac{|\l|^{m}-1}{|\l|-1}.$$ This however is impossible for $|\l|\ge 2$.

Now from the disjointness of $X_1=tX$ and $X_2=tvX$, and also $X_1\cup
X_2\subset X$, it now easily follows that $t, tv$ generate a free semigroup.
Any distinct words $u_1, u_2$ in  $t, tv$ without loss of generality begin on
the left with $t, tv$ respectively, so $u_1X\subset X_1, u_2X\subset X_2$ are
different. \endproof

\section{Abelian By Cyclic}

\begin{lemma} \label{two} Suppose that a finitely generated group $\G$ has a
normal free abelian subgroup $\V$, of rank $r$, with cyclic quotient generated
by $T$ having an eigenvalue of its characteristic polynomial of absolute value
at least 2, then $\G$  has uniform exponential growth. \end{lemma}

\proof We prove the proposition by induction on $r=rank_{\Z}(\V)$. Consider
$\V$ as a $Z[T]$ module with $T$ acting by conjugation. The eigenvalue
hypothesis of the statement implies that the  rank $r$ is at least 2. If $r=2$
then the invertible matrix $A$ has eigenvalues which are not integers. Consequently,
there are two generators of $\G$ $s,t$ which do not commute and $v=[s,t]\in
\V$; furthermore, one of them does not belong to $\V$, say $t$. Thus by Lemma
\ref{free}, $t, tv$ generate a free semigroup.

Suppose the rank $r$ is greater than 2, and inductively the lemma has been
shown for any rank less than $r$ which is at least 2. Given a generating set
$S$,  consider the collection $S_0$ of generators of $S$, which are in $\V$,
together with all commutators of generators.  Some generator, or its inverse,
say $t$, maps to a non-zero power of $T$ and so it also has an eigenvalue of
absolute value at least two.  The action of $t\in S$ on $\V$ satisfies its
characteristic polynomial of degree $r$, so in fact  the finite set of
conjugates of elements of $S_0$ by powers of $t^q$ for $0\le q< r$, yields
$S_1$ a generating set for the $t$ invariant subgroup $\V_1$.

If  $\V_1$ is of finite index in $\V$, the eigenvalue of absolute value at
least 2 will be supported on some cyclic submodule generated by one of the 
generators, since the least common multiple of all the annihilators of the
generators is the annihilator of the module, $\V_1$, which is the
characteristic polynomial of $T$.  Thus, by  Lemma \ref{free}, one of these
generators $v\in S_1$ of $\V_1$  generates a cyclic submodule, so that $t$ has
eigenvalue at least 2 on $Z[T]v\otimes_{\Z}\C$. The free semigroup generators
$t, tv$ have lengths bounded by $L=3+2r$, so $\beta(\G)\ge 2^{\frac{1}{L}}$.

If the rank of $\V_1$ is smaller than $r$ then $\V/\V_1$ has an action of $t$,
as well as $\V_1$; both groups are non-zero and have smaller rank than $\V$. 
To see that, suppose that $\V_1=\{0\}$ then all commutators of generators are
trivial and then $\G$ is abelian, which is of polynomial growth; this is
impossible by Lemma \ref{free}, since $\G$ must have exponential growth, and
so is not of polynomial growth. Thus $\V_1$ is non-zero.   Furthermore, the
torsion subgroup of $\V/\V_1$ is invariant under $t$ so that we may assume
$\V/\V_1$ is free abelian. On one or the other of these $\V_1$ or $\V/\V_1$,
the action of $t$ has an eigenvalue at least 2. To see this  we choose
rational bases for $\V_1$ and its complement, obtained by lifting a basis for
$\V/\V_1$. The characteristic polynomial of $t$ on $\V_1$ is a factor of that
on $\V$, with the complementary factor being given by the action on the
quotient group $\V/\V_1$.

If the eigenvalue at least 2 is supported on $\V_1$ then it must occur in the
exponent for this submodule, and hence is supported on one of the generators;
so we obtain the free semigroup as above with semigroup generators of length
bounded by a function of the rank and hence we have free semigroup generators
for $\G$ with generators of bounded length.

Otherwise, consider the subgroup $\G_1$ which is the inverse image of the
subgroup generated by $t$ in $\Z$; it contains $\V_1$ as a normal subgroup.
The group $\bar\G_1=\G_1/\V_1$ has a  (non-abelian) free semigroup by Lemma
\ref{free}, so by induction, $\bar\G_1$ is of uniform exponential growth;
hence also $\G$ has uniform exponential growth by Corollary $\ref{quo}$, since
$\bar\G_1$ is a homomorphic image of a subgroup of finite index. \endproof

The following important result of Kronecker is  well known.

\begin{proposition} \label{Kron} All the solutions of a monic integral
polynomial over $\C$ have absolute value at most 1 iff all the solutions are
roots of unity. \end{proposition} 
\medskip 

\begin{proposition}[{\it cf.}
\cite{R}] \label{abnil} Suppose that a finitely generated group $\G$ has a
normal free abelian subgroup $\V$, of rank $r$,  with nilpotent quotient,
$\N$, where every element acts so that all its eigenvalues are of absolute
value 1,  then  $\G$ is nilpotent by finite and hence of polynomial growth.
\end{proposition}

\proof Consider the matrix representation  afforded by the action of $\N$ on
the free abelian group $\V$. By the Lie-Kolchin-Mal'cev Theorem, \cite{K}, we
may put this integer representation in upper triangular form to obtain a
representation $\rho$. We construct the homomorphism
$\Delta:\rho(\N)\rightarrow (S^{1})^r$ mapping a matrix to its diagonal part.
Since each of the diagonal entries satisfies an integral polynomial, and all
of its conjugates are of absolute value 1 then by Proposition \ref{Kron} the
diagonal entries are roots of unity. Since the image of $\Delta$ is finitely
generated, it is of finite order. Thus $\rho(\N)$ has a subgroup $\rho(\U)$ of
finite index which is a group of unipotent matrices.  Hence, $\U/\V$ which is
a subgroup of finite index of $\N=\G/\V$ is also nilpotent. Moreover, we claim
that $\U$ is nilpotent and therefore it follows that $\N$ is nilpotent by
finite.

To see that $\U$ is nilpotent we continue an upper central series for $\V$;
denote the homomorphism $\phi: \U\rightarrow \U/\V$. Take a central series
$U_i=\phi^{-1}(Z_i)$ corresponding to $Z_i$ of $\U/\V$, so that $U_r=\U$.
Next, let $\{v_1,\cdots, v_s\}$ be the basis that makes $\rho(\U)$ unipotent.
Hence $\V$ is contained in the $\C$-span of $v_1,\cdots, v_s$.  Let $V_k=\V
\cap span_{\C}\{v_1,\cdots, v_k\}$. Then $V_s=\V$ and this gives a series
$V_i\subset V_{i+1}$ extending the series of $U_i$. Moreover, this is a
central series, because $V_{i+1}/V_{i}$ belongs to the center of $\U/V_{i}$.
To see this,  observe that $\rho(u)V_{i+1}\subset V_{i+1}$. Also,
$\rho(u)v_{i+1}=v_{i+1}\ mod\ span_{\C}\{v_1,\cdots, v_i\}$ for any element
$u\in \U$. Thus, for $v \in V_{i+1}$, then $ \rho(u)v-v\in V_i$. Hence $\U$ is
nilpotent and the polynomial growth for $\G$ follows from Wolf's Theorem.
\endproof

\begin{corollary} \label{abnil2} Suppose that a finitely generated group $\G$
has a normal free abelian subgroup $\V$, of rank $r$,  with nilpotent
quotient, $\N$. Then either $\G$ is nilpotent by finite or some element $t$ of
$\N$ acts on $\V$ so that not all its eigenvalues are of absolute value 1.  In
this later case  $\G$ is of exponential growth. \end{corollary} \proof In the
first case the result follows immediately from Proposition \ref{abnil}. In the
later case, the result follows from Lemma \ref{free}, using a large enough
power of $t$ or its inverse to get an eigenvalue of absolute value at least 2.
\endproof

\begin{theorem}[Abelian by Cyclic Alternative] \label{abelian} Suppose that a
finitely generated group $\G$ has a normal free abelian subgroup $\V$, of rank
$r$, with  cyclic quotient.  Then either $\G$ is of polynomial growth or is of
uniform exponential growth. \end{theorem}

\proof We suppose that the cyclic quotient is generated by $T$, and the module
action of $T$ is given by the matrix $A$. If all the eigenvalues of $A$ are
roots of unity then $\G$ is nilpotent by finite by Proposition \ref{abnil}. So
we shall assume otherwise.

The action of $T$ on the subgroup $\V$ satisfies its characteristic polynomial
which is of degree $r$. Since not all roots have absolute values bounded by 1,
we can assume that the matrix of $A^{n}$ has an eigenvalue of absolute value
at least  2, for some $n$, by Proposition \ref{Kron}.  Next, consider the
subgroup of $\G$ of finite index $\G_1=<\ \V,\ t^n >$, where $t$ maps to $T$.
It now follows by Proposition \ref{SW}  and Lemma \ref{two} that $\G$ has
uniform exponential growth. \endproof

\section{Polycyclic}

It is well-known that finitely generated polycyclic  group $\G$, \cite{K}, has
an embedding into $GL_m(Z)$. Whenever we have such an embedding and a normal
subgroup $\N$ we consider the action ( say $\rho$) of the group $\G$ (or possibly a 
quotient) as a group of automorphisms of the $\Z$-linear span of $\N$ in $M_m(\Z)$.
Consequently each automorphism has eigenvalues which belong to a ring of integers in a
number field of degree at most $m^2$ over the rationals. In fact, the eigenvalues of the
representation  $\rho(\gamma)$ are related (explicitly, $\lambda_i\lambda_j^{-1}$, $1\le
i,j\le m$ ) to the eigenvalues  ($\lambda_i$, $1\le i\le m$) of the matrix $\gamma\in
\G$.

Suppose that  $\G$ is a polycyclic group, which is an extension of a 
nilpotent normal subgroup $\N$ by a free abelian group $\A$ of finite rank.
Polycyclic groups always have a subgroup of finite index with this structure. 
Choose a
splitting of
$\A$ into
$\G$ and consider the action of
$\A$ on
$\G$ by inner automorphisms, using the matrix representation as above. This action
gives a sequence of actions on the terms of the lower central series of $\N$,
in which each element satisfies a polynomial of degree at most $m^2$, on each
term of the series. We refer to this as the action on $\A$ on $\N$.

The following important result is proved in \cite{W} or  \cite{G}, {\it appendix}.

\begin{lemma}
\label{rem}
With $\G$  above. For the action of 
 $\A$ ( or $\G$) on the lower
central series of
$\N$, there is an eigenvalue not of absolute value 1,  or else the group $\G$
is nilpotent by finite.
\end{lemma}

\begin{lemma} \label{polylemma} With $\G$ as above. Suppose that  the action of $\A$
on $\N$, every non-identity element, $t$ of $\A$, either $t$ or $t^{-1}$, has some
eigenvalue of absolute value at least 2, then $\G$ has uniform exponential
growth. \end{lemma} 

\proof We use induction on the Hirsh rank of $\N$,
$h=h(\N)$, the number of infinite-cyclic factors in a  cyclic  series for
$\N$. Consider a set of generators $S$ for $\G$. Let $\N_1$ be the normal
subgroup of $\N$, normally generated by the set $S_0$ of all the generators of
$S$ belonging to $\N$, together with all commutators of generators from $S$.
Since each element from $\G$ satisfies a polynomial of degree $m^2$, the
subgroup $\N_1$ is then generated by $S_0$, and the conjugates of this normal
set of generators by $t^q,\ 0\le q< m^2$,  using all the generators $t\in S$. Denote
this  set of generators of $\N_1$ by  $S_1$.It follows from the hypothesis, that $\G$ 
has exponential growth,
and hence it could not be virtually abelian; thus $\N_1$ is non-trivial.
If this nilpotent normal subgroup $\N_1$ has the same rank as $\N$ then
it will provide a free semigroup using the terms of the lower central series
and the generating set $S$ as we shall describe. The extension of $\N_1$ by $\A$ is of
finite index in $\G$ and as we show this is sufficient to finish the proof.
Let $\H=\N_1$. Consider the quotients of the lower central series
$\H_k/\H_{k+1}$ of $\H$. These are finitely generated abelian groups and
generated by the left normed commutators of weight $k$ in the generating set
$S_1$ for $\H$. Consequently, using Lemma \ref{rem}, some generator from $S$ acts on
the lower central series quotients so as to have an eigenvalue which is not a root of
unity. Since we have bounds on the lengths of a generating set of the lower
central series in terms of its weight and the lengths of elements of $S_1$
elements in terms of $S$ we will have a free semigroup, where the lengths of the
generators are bounded in terms of  $m$ and $h(\N)$.   The argument is
completed using Proposition \ref{bound} and Lemma \ref{free} which  guarantee that some
element of $t\in S\cup S^{-1}$ and $v$ either in  $S_1$ or a left normed commutator of
bounded length in elements of $S_1$  give a free semigroup. 

We continue with the induction argument.
If the rank of $\N_1$ is strictly smaller than the rank of $\N$,  and  $\A$ is
acting so that some element has an eigenvalue at least 2 in absolute value
then we finish the proof as in the previous case.

Otherwise, $\A$ acts so that every element has all eigenvalues which are roots
of unity. In this case the extension of $\N_1$ by $\A$ is virtually nilpotent by
Lemma \ref{rem}.  Now if both $\N/\N_1$ by $\A$ and $\N_1$ by $\A$ are nilpotent by
finite then we can construct the  representation of $\A$ on $\N$ which is the sum of
these two actions and hence has only eigenvalues which are roots of unity. In this
case $\G$ is nilpotent by finite, by Lemma \ref{rem}. This is
impossible since under the hypotheses of this lemma, $\G$  has exponential
growth (Lemma \ref{free}, Corollary \ref{abnil2} or see \cite{W}).  Hence, $\N/\N_1$
admits an action of
$\A$; the action here is such that some generator has an eigenvalue of absolute value
at least 2, and thus, by induction, the group extension of $\N/\N_1$ by $\A$ is of
uniform exponential growth; this group is a quotient of $\G$; hence, also $\G$ has
uniform exponential growth. \endproof

\begin{theorem}[Polycyclic Alternative] \label{poly} Suppose that  $\G$ is a
polycyclic by finite group, then either $\G$ has polynomial growth or has
uniform exponential growth. \end{theorem} \proof By using the
Lie-Kolchin-Mal'cev Theorem, \cite{K}, and passing to a subgroup of finite
index, $\G$ can be realized as a torsion free group of upper triangular
matrices, which has commutator subgroup, $\N$,  which is nilpotent,  and
represented by unipotent matrices. We may  assume that $\G$ is also poly-$\Z$
by passing again to a subgroup of finite index.  The subgroup $\N$  is
finitely generated, since polycyclic groups satisfy the maximum condition on
subgroups (Noetherian).
We consider the inverse image, $\N_1$, in $\G$, of the set of
elements of $\G$ where the elements act with all eigenvalues, roots of unity; this
subgroup,
$\N_1$, is nilpotent by finite by Lemma \ref{rem}; it contains the
commutator subgroup of $\G$; it is normal in $\G$ with quotient which is the
finitely generated abelian group, $\A$.
The inverse image $\G_1$ in $\G$ of the maximal torsion-free summand $\B$ in
$\A$ is of finite index in $\G$. If the torsion free rank of $\B$ is zero,
then of course $\G$ is nilpotent by finite. If the torsion free rank of $\B$
is nonzero, It suffices by Proposition \ref{SW} to show $\G_1$ is of uniform
exponential growth.

We then may assume then that $\G_1$ is nilpotent by free abelian, where every
non-identity element of the non-trivial free abelian group $\A$ acts with an
eigenvalue which is not a root of unity. Choose $n$, by Proposition \ref{Kron}
so that the eigenvalues of every non-trivial element of $\A$ or its inverse
has $n^{th}$ power which has an eigenvalue which has absolute value greater
than or equal to 2. The $n^{th}$ power of all elements in $\A$ forms a
subgroup of finite index in $\A$ and its inverse image is a subgroup of finite
index in $\G_1$ satisfying the hypotheses of Lemma \ref{polylemma}. Hence it
is of uniform exponential growth.

\section{Solvable}

We denote the terms of the derived series of $\G $  by $\G^{(k)}$. When
$\G^{(k-1)}/\G^{(k)}$ are finitely generated for all $k\le m$ then
$\G/\G^{(m)}$ is polycyclic.  The following is an immediate corollary from
Theorem \ref{poly}.

\begin{corollary} Suppose that  $\G$ is a finitely generated  group, for which
some quotient, $\G/\G^{(m)}$, is polycyclic but not nilpotent by finite then
$\G$ has uniform exponential growth. \end{corollary}

\section*{Concluding Remarks}

The results of this paper were announced at the
Stallingsfest Group Theory meeting in May 2000 at M.S.R.I. in Berkeley. Six weeks
later, I learned that Denis Osin (preprint: `The entropy of solvable
groups') had obtained the same results independently and at essentially the
same time; moreover he  has also extended the polycyclic case to prove that any
finitely generated solvable group of exponential growth has uniform exponential growth.
Subsequently, another proof of Osin's result has been given by J. S. Wilson.

\end{document}